\documentclass[12pt]{article}

\usepackage{graphicx}
\usepackage{bm}
\PassOptionsToPackage{hyphens}{url}\usepackage{hyperref}
\usepackage{mathtools}
\usepackage{subcaption}
\usepackage{tikz}
\usepackage{tabto}
\usepackage{arydshln}

\input epsf
\usepackage{amsmath}
\usepackage{amssymb}
\usepackage{amsthm} 
\usepackage[margin=2cm]{geometry}

\usepackage{algorithm2e}

\newtheorem{theorem}{Theorem}[section]

\newtheorem{remark}[theorem]{Remark}

\usepackage{moreverb}

\newcounter{gagcomment}

\newcommand{\R}{\mathbb{R}}
\newcommand*{\Id}{\operatorname{Id}}
\newcommand*{\X}{{\boldsymbol{X}}}
\DeclareMathOperator*{\argmin}{arg\,min}
\def\x{{\bf{x}}}

\newcommand*{\bchi}{{\boldsymbol{\chi}}}
\newcommand*{\bchis}{{\hat \bchi^{(s)}}}

\title{SINDy for delay-differential equations: application to model bacterial zinc response 
}

\author{
Antoine Sandoz\thanks{Section of Mathematics and Microbiology Unit, 
Department of Plant Sciences,
University of Geneva, Switzerland; {\rm{Antoine.Sandoz@unige.ch}}} 
\and Verena Ducret\thanks{Microbiology Unit, 
Department of Plant Sciences,
University of Geneva, Switzerland; {\rm{Verena.Ducret@unige.ch}}}
\and Georg A. Gottwald\thanks{School of Mathematics and Statistics, University of Sydney, Sydney, Australia; {\rm{georg.gottwald@sydney.edu.au}}} 
\and Gilles Vilmart\thanks{Section of Mathematics, University of Geneva, Switzerland; {\rm{Gilles.Vilmart@unige.ch}}}
\and Karl Perron \thanks{Microbiology Unit, 
Department of Plant Sciences, and Section of Pharmaceutical Sciences,
University of Geneva, Switzerland; {\rm{Karl.Perron@unige.ch}}}
}

\date{\today}

\begin{document}
\maketitle
\abstract{We extend the data-driven method of Sparse Identification of Nonlinear Dynamics (SINDy) developed by {\em{Brunton et al}, Proc. Natl. Acad. Sci USA 113 (2016)} to the case of delay differential equations (DDEs). This is achieved in a bilevel optimization procedure by first applying SINDy for fixed delay and then subsequently optimizing the error of the reconstructed SINDy model over delay times. We test the SINDy-delay method on a noisy short data set from a toy delay differential equation and show excellent agreement. We then apply the method to experimental data of gene expressions in the bacterium {\em{Pseudomonas aeruginosa}} subject to the influence of zinc. The derived SINDy model suggests that the increase of zinc concentration mainly affects the time delay and not the strengths of the interactions between the different agents controlling the zinc export mechanism.}

\vspace{1cm}

\noindent
{\bf{keywords}}: data-driven modelling; SINDy; delay-differential equations; {\em{Pseudomonas aeruginosa}}, zinc homeostasis\\

\section{Introduction}
Our ability to understand, forecast and control dynamical systems depends crucially on our knowledge of its underlying equations. Recently data-driven methods to uncover underlying equations 
have been proposed to uncover unknown dynamics and to increase our forecast capabilities 
\cite{BruntonKutz,ChampionEtAl18,SchmidtLipson09,GuimeraEtAl20,UdrescuTegmark20,BiggioEtAl21}. 
A particularly attractive and easy-to-implement method, proposed by  Brunton et al. \cite{BruntonEtAl16}, is {\em{Sparse Identification of Nonlinear Dynamics}}, or  SINDy for short. The problem which is addressed by SINDy is the following: Given observations $\x_n\in \R^d$ sampled at (not necessarily equidistant) times $t_n$ which were generated by a dynamical system of the form 
\begin{align*}
\dot \x = \mathcal{F}(\x),
\end{align*}
where the dot signifies the time derivative,
find an approximation to this dynamical system using only the data. SINDy approaches this question by assuming that the vector field $\mathcal{F}(\x)$ lies in the span of a given  (potentially very large) library of functions such as simple polynomials. This reduces the problem to linear regression on a library of nonlinear functions. In line with the parsimony principle, SINDy imposes a sparsity constraint leading to a sparse approximation of $\mathcal{F}(\x)$ in terms of the members of the library functions. If the system is only partially observed, the method of SINDy can be extended to the reconstructed phase-space using Takens' embedding theorem and it is then known as the Hankel alternative view of Koopman (HAVOK) analysis \cite{BruntonEtAl17}. SINDy has been successfully applied to systems appearing in a wide range of scientific disciplines, including fluid dynamics, plasma physics and nonlinear optics \cite{LoiseauBrunton18,DamEtAl17,SorokinaEtAl16}.

Many dynamical systems involve a delayed feedback response and are modelled by delay-differential equations (DDEs). Examples range from the natural world to engineering with applications in, for example, population dynamics \cite{Cushing77,Gourley00}, biological regulatory systems \cite{GlassEtAl21}, cardiac dynamics \cite{GottwaldKramer06,Gottwald08}, climate dynamics \cite{SuarezSchopf88,KeaneEtAl17}, mechanical vibration \cite{WangEtAl19} and in optical systems \cite{TerrienEtAl21}, to name just a few. 
In this paper, we extend the framework of SINDy to dynamical systems which are described by DDEs and where in addition to the sparse subset of the library of nonlinear functions and their associated coefficients, the delay has to be determined as a parameter. We achieve this by employing a bilevel optimization in which, for a fixed specified delay time, the error in reproducing the observations made by each approximate SINDy model is minimized. Particular emphasis will be given to deal with noisy data.
We shall first test the proposed SINDy-delay methodology to a one-dimensional toy model with artificially noisy data before considering a challenging 
problem with biological data of gene expressions in the bacterium {\em{Pseudomonas aeruginosa}} subject to the influence of zinc.

\textit{P. aeruginosa} is an opportunistic pathogen capable of causing acute infections in hospitals, in particular in immunocompromised patients, in cystic fibrosis patients and in severe burn victims \cite{Kerr09,Jurado-Martin21}.
Therefore, it belongs to the Priority 1 category for research into antibiotic resistance as determined by the world health organization \cite{WHO2018Lancet}. \textit{P. aeruginosa} has a large genome, coding for 5570 open
reading frames, of which 72 are involved in predicting so called two component systems (TCS) \cite{stover_complete_2000}. TCSs are crucial biological building blocks. They are 
composed of a sensor protein which in response to a stimulus, activates a cognate transcriptional
regulator by phosphorylation, allowing for a rapid adaptation to environmental changes. \textit{P. aeruginosa} has one of the highest numbers of putative TCSs among bacteria,
contributing to the ubiquity of this micro-organism \cite{alm2006}. For instance, the CzcRS TCS promotes resistance to high concentrations as well as to large fluctuations of the concentration of trace metals such as zinc. This is advantageous for the bacteria in an infectious context \cite{Perron04,Dieppois12} since to counter the multiplication of bacteria, the host uses
nutritional immunity strategies via scavenging essential nutrients including zinc,
iron and manganese \cite{kehl-fie_nutritional_2010, capdevila2016,lonergan2019}. 
Conversely, during phagocytosis, macrophages deliver a toxic amount of
zinc and copper into the phagolysosome, leading to the death of the invader organism \cite{Djoko15,stafford2013,gao2018}. Thus, the
success of an infection depends largely on the capacity of a pathogen to survive in zinc deficient
as well as zinc excess environments and to switch from one to the other of these extreme
conditions. \textit{P. aeruginosa} has a whole arsenal of the most effective systems for regulating the
entry and exit of the metal. Moreover, zinc was shown to exacerbate the bacterium pathogenicity,
enhancing the virulence factor production and rendering this micro-organism more resistant to
antibiotics, especially those belonging to the carbapenem family, a last resort anti-pseudomonas
class of compounds \cite{Perron04,Dieppois12}. In order to better understand the \textit{P. aeruginosa} zinc homeostasis, we derive a mathematical delay differential equation model focusing on the dynamic of two main zinc export machineries. To address this challenging question, we use the proposed SINDy-delay methodology applied to experimental data. 

The paper is organised as follows. In Section~\ref{sec:SINDy} we introduce an extension of SINDy to find parsimonious models for delay-differential equations. In Section~\ref{sec:ENSO} we illustrate the effectiveness of our method in the context of a known toy model DDE; in particular, we show that the DDE is recovered well even for short noise-contaminated observations. In Section~\ref{sec:bio} we then apply our method to experimental data of gene expression data of the bacterium \textit{P. aeruginosa} under various concentrations of zinc. We conclude in Section~\ref{sec:conclusion} and discuss biological implications of the discovered DDE describing the bacterium's zinc regulation system.


\section{Sparse Identification of Nonlinear Dynamics with delay for noisy data (SINDy-delay)}
\label{sec:SINDy}
Consider a $d$-dimensional dynamical system with delay time $\tau$,
\begin{align} \label{eq0}
\dot \x = \mathcal{F}(\x(t),\x(t-\tau)),
\end{align}
where $\x(t)\in\R^d$, which is probed at times $t_n$, $n=\ldots, -1, 0,1,2,3\ldots,$ by observations
\begin{align*}
\boldsymbol{\chi}_n = \x_n +\boldsymbol{ \Gamma} \boldsymbol{\eta}_n,
\end{align*}
with measurement error covariance matrix  $\Gamma^2 \in \R^{d\times d}$ and independent normally distributed noise ${\boldsymbol{\eta}}_n\sim {\mathcal{N}}(0,\Id)$. For simplicity we assume here throughout $\boldsymbol{ \Gamma}  = \gamma \Id$. The aim is to find a  parsimonious approximation of the vector field $\mathcal{F}(\x(t),\x(t-\tau))$ as a linear combination of nonlinear functions selected from a library $\mathcal{R}$ of cardinality $N_{\mathcal{R}}$. In particular, the $k$th component is expressed as a linear combination of all the library functions $\theta_{j} \in \mathcal{R}$,  $ j=1,\dots,N_{\mathcal{R}}$, of the form
\begin{align} 
\mathcal{F}_k(\x(t),\x(t-\tau))=\sum_{j=1}^{N_{\mathcal{R}}} \xi_{k}^{j} \theta_{j}(\x(t),\x(t-\tau)) + \epsilon_k(\x(t),\x(t-\tau)),
\label{e.eps}
\end{align}
for all components $k=1,\ldots,d$. Simple nonlinear regression would amount to determining the coefficients $\xi_{k}^{j}$ using the method of least-squares to minimize the mismatch $\epsilon_k$. In SINDy rather, a sparsity constraint is invoked, seeking a parsimonious model with as many of the coefficients $\xi_{k}^{j}$ being zero while still ensuring fidelity of the approximation (\ref{e.eps}) with respect to the data.  

To describe how SINDy finds such an approximation, let us assume for the moment that observations are taken at equidistant times $t_n=n\Delta t$ with constant sampling time $\Delta t$. To account for the delay we form the observation vector
\begin{align*} 
\hat \bchi_n^{(s)} = 
\begin{pmatrix*}[l]
\bchi_n\\
\bchi_{n-s}
\end{pmatrix*}
\in \R^{2d},
\end{align*}
for $n=1,\dots,N$, where the positive integer $s$  is related to a delay time $\tau=s\Delta t$. Following the exposition in Brunton et al. \cite{BruntonEtAl17} and Brunton and Kutz \cite{BruntonKutz}, we collect the observation vectors $\hat\bchi^{(s)}_n$ in a data matrix
\begin{align*}
\boldsymbol{X}^T 
= \begin{pmatrix*}[l]
\bchis_1 & \bchis_2&\dots&\bchis_{N}
\end{pmatrix*}
\in \R^{2d \times N}.
\end{align*}
Similarly we define the matrix consisting of derivatives of the observations $\bchi_n$ at the observation times 
 \begin{align}\label{eq:Xdot}
\boldsymbol{\dot X}^T  
= \begin{pmatrix*}[l]
{\dot { \bchi}}_1 & {\dot {\bchi}}_2&\dots&{\dot {\bchi}}_{N}
\end{pmatrix*}
\in \R^{d \times N}.
\end{align}
Note that for the derivatives we only consider the time derivatives for $\bchi_n$ and not for $\bchi_{n-s}$ which would be redundant information. Typically one does not have access to the actual derivatives $\dot \bchi_n$ but only to the variables $\bchi_n$ themselves. For noise-free finely sampled observations with $\Delta t\ll 1$ finite differencing can be employed to approximate the derivatives. For noisy observations, however, estimating derivatives via finite-differencing leads to an amplification of the noise. Denoising methods such as the total-variation regularized method are required \cite{Chartrand11}. Here we propose to use simple polynomial regression for denoising as discussed in the following remark. 

\begin{remark} \label{rem:noise} (Denoising procedure for computing the derivative matrix \eqref{eq:Xdot})
For computing each $\dot \bchi_n$ , $n=1\ldots, N$ in  \eqref{eq:Xdot}, we use polynomial regression. We define for the corresponding observation time $t_n$ a temporal window $[t_n-\delta,t_n+\delta]$ with $\delta = r\Delta t$ and fit a $3$rd order polynomial through the $2r-1$ observations $\bchi_{n-r},\bchi_{n-r+1},\ldots,\bchi_{n+r}$ lying within this time window. Choosing a sufficiently large temporal window containing more data points compared to the regression polynomial degree allows for noise reduction. 

The derivatives $\dot \bchi_n$ can then be analytically determined from the fitted polynomials at each time $t_n$. This denoising procedure can easily be adapted to handle non-equidistantly sampled observations which may be the situation in experimental data (including the case of determining delayed data at  $t_n-\tau$ which may not have been directly observed). Our denoising procedure is closely related to the Savitsky-Golay filter \cite{SavitzkyGolay64} and denoising by splines \cite{Wahba75}; for a comparison of various denoising strategies see \cite{VanBreugelEtAl20}.
\end{remark}

At the heart of SINDy lies the choice of a suitably large library $\mathcal{R}$. A natural choice is the set of monomials in $x_k(t),x_k(t-\tau),k=1,\ldots,d$ up to a fixed degree $M$, $\mathcal{R}=\{1,x_1(t),x_2(t),x_1(t-\tau),x_2(t-\tau)\ldots \}$ with cardinality $N_{\mathcal{R}} = \binom{2d+M}{M}$. Given a library $\mathcal{R}$, the associated library matrix $\boldsymbol{\Theta} (\X)\in \R^{N \times  N_{\mathcal{R}}}$ is constructed from the data by evaluating all functions $\theta_{j}(\x(t),\x(t-\tau))$ of the library $\mathcal{R}$ at the observation times $t=t_1,\ldots, t_N$. When considering the library consisting of monomials of up to order $M$, the library matrix becomes
 \begin{align*}
\boldsymbol{\Theta}(\X)
= \begin{pmatrix}
\boldsymbol{1} & \X & \X^2 & \X^3 & \ldots &  \X^M
\end{pmatrix}
,
\end{align*}
where the matrices $\X^m \in \R^{ N \times \binom{2d+m-1}{m}}$ consist of rows whose coefficients include all possible monomials of degree $m$ between the $d$-dimensional variables $\bchi_n$ and $\bchi_{n-s}$. For simplicity, we will later in the numerical experiments exclude any products between $\bchi_n$ and $\bchi_{n-s}$. This reduces the number of columns of each $\X^m$ to $2\binom{d+m-1}{m}$ and the overall number of columns of $\boldsymbol{\Theta}(\X)$ to $N_{\mathcal{R}} =  2\binom{d+M}{M}$.

In SINDy the minimization of the error $\epsilon_k$ made by the approximation (\ref{e.eps}) is achieved by an $\ell_1$-regularized regression problem. Defining first the $\ell_2$-cost function
\begin{align}
C(\Xi) = \sum_{k=1}^d \| {\dot{X}}_k - \boldsymbol{\Theta}(\X) \xi_k\|_2^2,
\label{e.C}
\end{align}
where $\dot X_k \in \R^N$ denotes the $k$th column of $\dot \X$ and $\Xi=\{\xi_k\}_{k=1,\ldots,d}$ is the coefficient matrix consisting of column vectors $\xi_k\in \R^{N_{\mathcal{R}}}$ which denote the coefficients associated with the library functions for the $k$th component of the state variable (cf. (\ref{e.eps})). To promote sparsity of the coefficients the cost function is minimized under an $\ell_1$-sparsity constraint according to
\begin{align}
\boldsymbol{\xi_k} = \argmin_{\xi_k\in \R^{N_{\mathcal{R}}}} \| {\dot{\X}}_k - \boldsymbol{\Theta}( \X) \xi_k\|_2+ \lambda \|\xi_k\|_1,
\label{e.cost}
\end{align} 
where the regularization parameter $\lambda$ controls the sparsity. Rather than using a sequential thresholded least-squares algorithm to approximate the solution of the optimisation problem (\ref{e.cost}), as suggested in \cite{BruntonEtAl16}, we promote here sparsity by the following sequential procedure. Define $\xi^q_k$ to be the coefficient of the $q$th library function $\theta_q$ which is associated with the $k$th component of the vector field. For each $q=1,\dots,N_{\mathcal{R}}$  calculate the least square solution $\xi_k^{q} \in \R^{N_\mathcal{R}}$ corresponding to the minimization of the cost function $C(\Xi)$ with the hard sparsity constraint
\begin{align*}
\xi_{k}^q=0.
\end{align*}
For each of the $N_{\mathcal{R}}$ solutions $\xi_k^{q}$ record the associated minimized cost $C(\Xi)$, and select the value 
\begin{equation} \label{eq:qstar}
q^*=\argmin_{q=1,\ldots,N_{\mathcal{R}}} \min_{\xi_k} \{C(\Xi)\ ; \ \xi_{k}^q=0 \}
\end{equation}
corresponding to the hard sparsity constraint $\xi_{k}^{q^*}=0$ which leads to the smallest increase in the minimum of the cost $C(\Xi)$.  We then set $\xi_{k}^{q^*}=0$, i.e. excluding $\theta_{q^\star}$ from the library $\mathcal{R}$ for the $k$th state variable. Algorithmically this amounts to deleting the $q$th column of $\boldsymbol{\Theta}( \X)$ when seeking solutions of (\ref{e.cost}). 
This process of eliminating coefficients $\xi_k^q$ is then repeated for the remaining library functions in $\mathcal{R}$ (and the corresponding columns of $\boldsymbol{\Theta}( \X)$) until a significantly large change of the cost $C(\Xi)$ has been accrued, suggesting that removing any of the remaining functions will lead to a strong increase of the cost function, thereby deteriorating the accuracy of the SINDy model. 

\begin{remark} \label{rem:sparsity}
(Promoting sparsity to approximate solutions of (\ref{e.cost})) 
Promoting sparsity by envoking \eqref{eq:qstar} avoids having to set a cutoff value $p$ such that coefficients with $|\xi_k^j|\leq p$ are removed as proposed in Brunton et al. \cite{BruntonEtAl16}. Instead the degree of sparsity is visually determined by plotting the cost function for an increasing number of removals. This does not require the data $X$ to be normalized in a pre-processing step and can be applied to situations in which variables may exhibit widely varying ranges. We shall encounter such a situation for the experimental data in Section~\ref{sec:bio}. 
\end{remark}

The above procedure is applied to each of the components $k=1,\ldots,d$ with each component having their separate subset of eigenfunctions selected. Collecting the typically sparse output vectors $\xi_k^*$, $k=1,\ldots d$ in the matrix $\Xi^*=(\xi_1^*,\ldots \xi_d^*)$, the approximate SINDy delay differential equation model for arbitrary fixed delay time $\tau_s=s\Delta t$ is given by
\begin{equation} 
\label{e.SINDy}
\dot x_k(t;\tau_s) = \sum_{j=1}^{N_{\mathcal{R}}} (\xi^j_k)^* \theta_{j}(\x(t),\x(t-\tau_s)),\quad k=1,\ldots,d. 
\end{equation}
Up to here this is standard SINDy, as described in Brunton et al. \cite{BruntonEtAl16}, except for the proposed alternative method of denoising with local polynomial regressions of the data points described in Remark \ref{rem:noise}, and for the modified algorithm to approximate solutions to the optimization problem (\ref{e.cost}) described in Remark \ref{rem:sparsity}. 

\begin{algorithm}[tb] 
\hrule
\smallskip
 \textbf{Input:} {Observational data $\chi_n,n=\ldots,-1,0,1,2,\ldots$, and nonlinear function library set $\mathcal{R}=\{\theta_j\ ;\ j=1,\dots, N_{\mathcal{R}} \}$}.\\
  \textbf{Output:} {Delay time $\tau^*$ and coefficient matrix $\Xi^*$ determining the delay differential equation model \eqref{eq_S}}.\\[2mm]
 Compute the derivative matrix $\boldsymbol{\dot X}$ in \eqref{eq:Xdot}, with denoising procedure (see Remark~\ref{rem:noise})\;
 \For{all delay time $\tau_s=s\Delta t$, $s\in\{0,1,2,\ldots\}$ }{
 Compute the data matrix $\boldsymbol{X}$ and the associated library matrix $\Theta(\boldsymbol{X})$.\\
 \For{$k\in\{1,\ldots,d\}$ }{
  Set the list $Q\subset\{1,\ldots,N_{\mathcal{R}}\}$ of indices of vanishing coefficients $\xi_k^j=0$ to achieve the sparsity constraint of the SINDy methodology to $Q=\emptyset$\;
	\While{$C=\min_{\xi_k \in \R^{N_\mathcal{R}}} \{ \| {\dot{X}}_k - \boldsymbol{\Theta}(\X) \xi_k\|_2 \ ;\ \xi_k^j =0\mbox{ for all }j\in Q\}$ does not increase significantly (see Remark~\ref{rem:sparsity})}{
	Compute 
	$\displaystyle
	(q^*,\xi_k^*)=\argmin_{q\in\{1,\ldots,N_{\mathcal{R}}\}\setminus Q,\ \xi_k \in \R^{N_{\mathcal{R}}}}\{\| {\dot{X}}_k - \boldsymbol{\Theta}(\X) \xi_k\|_2\ ; \ \xi_{k}^j=0 \mbox{ for all } j\in Q \cup \{q\}\}
	$.\\
	Set $Q=Q\cup \{q^*\}$.
	}
	}
	Keep $\Xi^*=\{\xi_1^*,\ldots,\xi_d^*\}$ and the corresponding error $\mathcal{E}(\tau_s)$ in \eqref{e.ell2}.
 }
 Save the optimal delay time $\tau^*$ given in \eqref{e.cost2}, and the corresponding coefficient matrix $\Xi^*$.
\smallskip
\hrule
\smallskip
 \caption{SINDy algorithm for dynamical systems with temporal delay.
\label{algo1}}
\end{algorithm}

To account for a delay we extended the nonlinear library $\{\theta_k\}_{k=1,\dots,{N_{\mathcal{R}}}}$ to include delay terms $\x(t-\tau)$, which fits in the standard SINDy methodology for fixed delay time parameter $\tau_s=s\Delta t$. To estimate the delay time $\tau=s\Delta t$ of the dynamical system \eqref{e.SINDy} which best matches the data $\boldsymbol{\chi}_n$ an additional optimization procedure is employed: Consider a range of delay times $\tau_s = s\Delta t$ with integer sequence $s\in \{0,1,2,\ldots,\}$. For each $\tau_s$ we perform the above procedure to obtain the SINDy model (\ref{e.SINDy}). We then compute the reconstruction error $\mathcal{E}(\tau_s)$ for fixed delay time parameter $\tau_s$ as the $\ell_2$-error between the solution of the SINDy model (\ref{e.SINDy}) and observations $\bchi_n$,
\begin{align}
\mathcal{E}(\tau_s)=\frac{1}{Z}\sum_{n=1}^N \|\x(t_n; \tau_s) - \bchi_n\|^2.
\label{e.ell2}
\end{align}
We set the normalization constant to $Z=\sum_{j=0}^N \|\bchi_j\|^2$. Note that we define the error $\mathcal{E}(\tau_s)$ here in terms of the trajectories $\x(t)$ rather than via the derivatives as in the cost function (\ref{e.C}) used in the SINDy core. This proves to be advantageous as trajectories are less affected by the noise than their derivatives. The optimal delay time is estimated as the solution of  
\begin{align}
\tau^\star = \argmin_{\tau_s} \mathcal{E}(\tau_s).
\label{e.cost2}
\end{align} 
This finally yields the SINDy-delay differential equation model,
\begin{align} 
\label{eq_S}
\dot \x(t)^T = \boldsymbol{\Theta}( \x(t)^T,\x(t-\tau^*)^T)\Xi^*.
\end{align} 
We remark that the devising strategy proposed in Remark~\ref{rem:noise} allows for non-uniformly sampled data, provided the sampling times are not too far apart. The required values of unobserved $x(t-\tau)$ can be evaluated by the proposed polynomial regression.
We summarize the extension of the SINDy methodology to systems involving temporal delays in Algorithm \ref{algo1}. 
In the next Section we show how the SINDy-delay method performs for artificial data obtained from a simple one-dimensional DDE.


\section{Application to a toy model}
\label{sec:ENSO}

\begin{table}[btp]
\begin{center}
\small
\begin{tabular}{ccccccccccc}
&  $N$& $\Delta t$& noise $\gamma$ & observations &  $\tau^\star$ & $x(t)$ & $x(t-\tau)$  & $x^3(t)$ & error $\mathcal{E}(\tau^{\tau^\star})$\\
\hline
&--& -- & --& -- &                             $7$ & $1$ & $-0.75$ &  $-1$ & -- \\
\hline
(a) & $4000$ & $0.025$ & $0$     & $x,\dot{x}$ &  $7.00$  &$1.00$	&$-0.75$	&$-1.00$&		$1.36e-03$\\
(b) & $4000$ & $0.025$ & $0$     & $x$ & $7.00$ &$0.99$	&$-0.75$		&$-0.99$		&$1.80e-03$\\
(c) & $200$ & $0.25$  & $0.02$ & $x$ & $7.00$ &$0.84$	&$-0.71$		&$-0.88$		&$2.56e-02$\\
(d) & $4000$ & $0.025$ & $0.02$ & $x$ & $7.025$ &$0.92$	&$-0.74$		&$-0.94$		&$1.94e-02$\\
\hline
\end{tabular}
\end{center}
\caption{Results of the SINDy-delay method (Algorithm \ref{algo1}) applied to data obtained from the toy model (\ref{e.ENSO}). The first row denotes the true delay time and coefficients of the DDE (\ref{e.ENSO}) used to generate the observations. Rows (a)--(d) present results for the different scenarios described in the main text. We present results the estimated delay times $\tau$, the coefficients as well as the associated reconstruction error $\mathcal{E}(\tau)$ (cf \eqref{e.ell2}) between the data and the SINDy differential equation model \eqref{eq_S} for varying data length $N$, sampling times $\Delta t$, noise levels $\gamma$. The columns for the monomial terms $x(t), x(t-\tau),x(t)^3$ display the estimated corresponding coefficients $\xi_{1}^j$ (coefficients for the remaining monomials were estimated to be $0$ in all experiments). Experiments in which only $x$ is observed required polynomial regression as described in Remark~\ref{rem:noise}.}
\label{tab:para}
\end{table}

To illustrate how the SINDy-delay method is able to determine an underlying DDE together with the delay time parameter from noisy observations, we consider the following one-dimensional DDE,
\begin{align}
\dot x(t) = x(t)-x(t)^3 - \alpha x(t-\tau).
\label{e.ENSO}
\end{align}
This DDE was introduced as a toy model in the context of climate science to describe, for example,  the El Ni\~no -- Southern Oscillation (ENSO) phenomenon where $x(t)$ denotes a sea-surface temperature anomaly at time $t$ \cite{SuarezSchopf88}. We choose in the following as parameter value $\alpha=0.75$ and a delay time of $\tau=7$. The initial solution on the time interval $[-\tau,0]$ is chosen to be the stable periodic solution to \eqref{e.ENSO} and with $x(0)=1$. 
For the set of nonlinear library functions $\mathcal{R}$, we consider all monomials up to cubic degree, $\mathcal{R}=\{1,x(t),x(t)^2,x(t)^3,x(t-\tau),x(t-\tau)^2,x(t-\tau)^3\}$, excluding products of $x(t)$ and $x(t-\tau)$. We simulate the DDE (\ref{e.ENSO}) using the Matlab dde23 integrator with absolute and relative tolerances of $10^{-8}$ to produce time series of $N$ observations sampled at equidistant times with sampling time $\Delta t$ \cite{MATLAB:2019}.

We present results for several scenarios with increasing difficulty. In particular, we investigate how the accuracy of the method depends on the amount of data available as well as on the level of noise. In the following we restrict the delays $\tau_s$ to $\tau_s = s \Delta t$ for $s=1,\dots,8.5/\Delta t$ so that we sample from the interval $[0, 8.5]$. 

\begin{figure}[tbp]
\centering
\begin{subfigure}[t]{0.47\textwidth}
			\centering
			\includegraphics[width=1\linewidth]{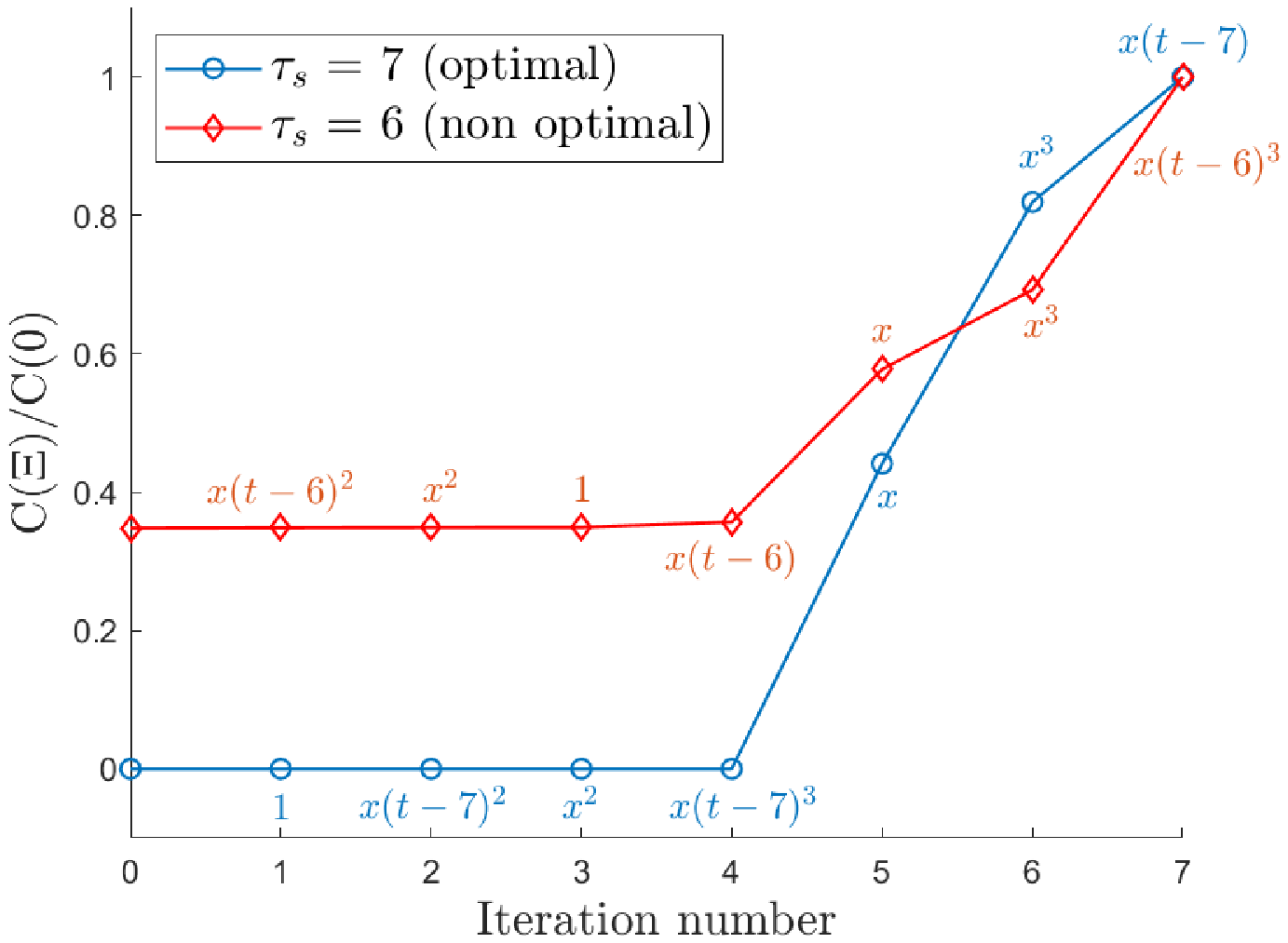}
			\caption{noiseless observations of $x$ and $\dot{x}$ ($N=4\,000$).}
		\end{subfigure}
\begin{subfigure}[t]{0.47\textwidth}
			\centering
			\includegraphics[width=1\linewidth]{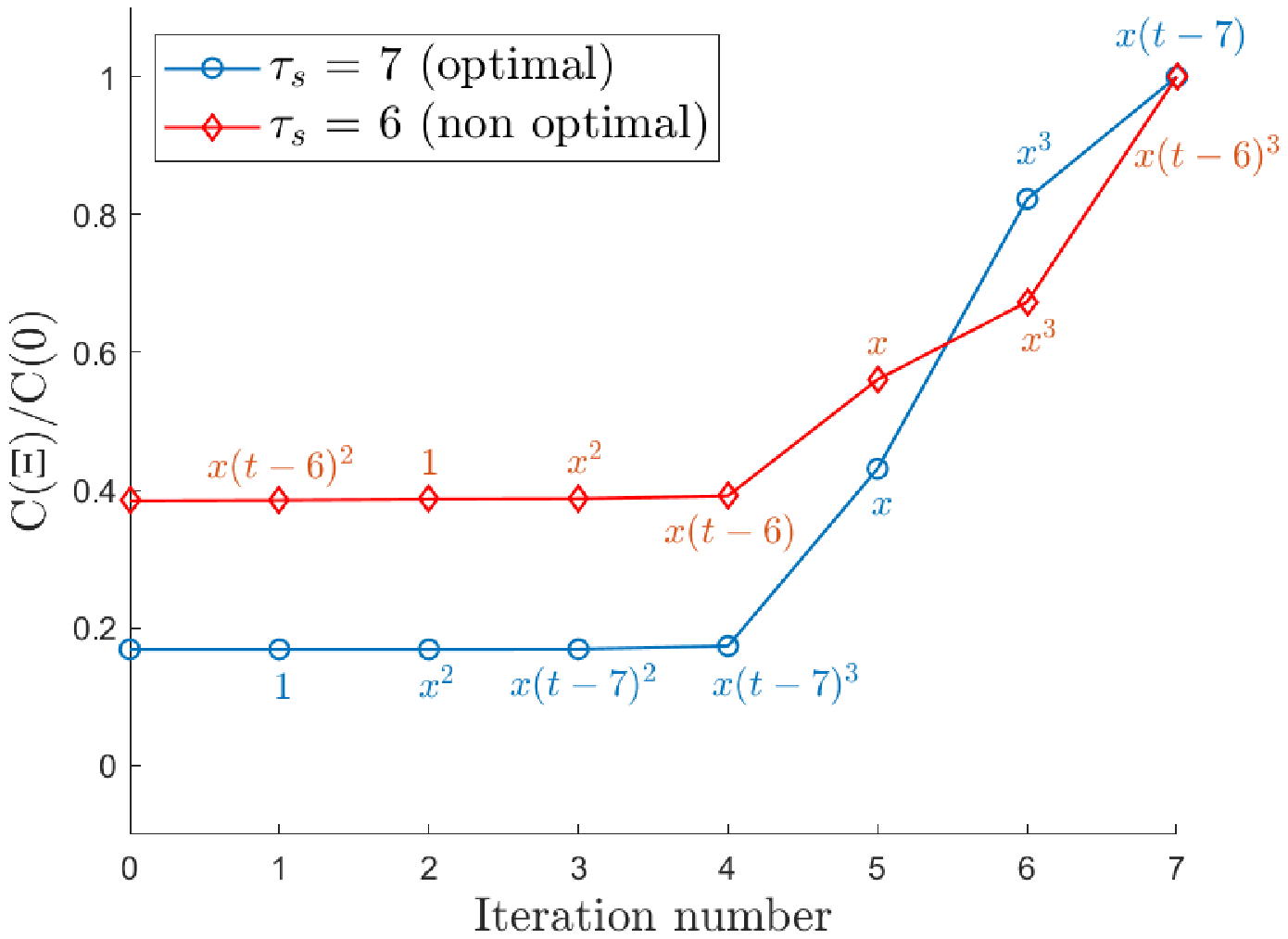}
			\caption{noisy observations of $x$ ($N=200$, $\gamma=0.02$).}
		\end{subfigure}
\caption{Cost function $C(\Xi)$ (nomalized by $C(0)$) against removed monomials for fixed delay. Results are shown for the true delay time $\tau_s=7$ (open circles, online blue) and for the non-optimal delay time $\tau_s=6$ (diamonds, online red). The error increases after iteration 4 of the SINDy algorithm, as soon as the terms $x,x^3,x(t-\tau)$ actually present in the underlying DDE model start to be removed from the library.}
\label{fig:remove}
\end{figure}

\begin{figure}[tbp]
\centering
\includegraphics[width=0.5\columnwidth]{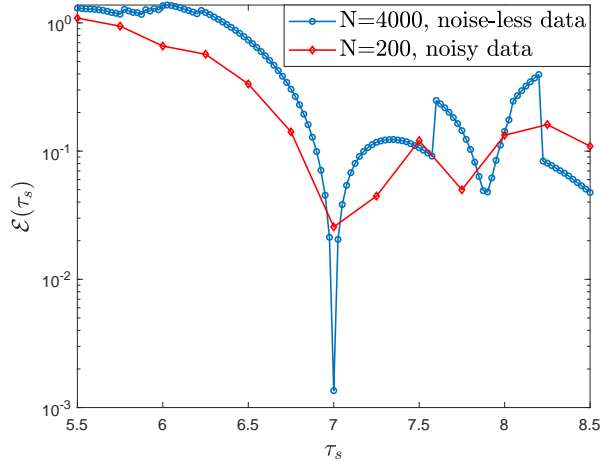}
\caption{Reconstruction error $\mathcal{E}(\tau_s)$ as a function of the delay time $\tau_s$, showing a clear minimum at the true delay time $\tau=7$ for the ideal case with $N=4\,000$ noiseless observations of $x$ and $\dot{x}$ (open circles; online blue) and for the case of $N=200$ noisy observations of $x$ with noise level $\gamma=0.02$ (diamonds; online red).}
\label{fig:L2error}
\end{figure}

We begin with the ideal situation of noiseless observations of both the state $\x$ and the derivative $\dot\x$. We consider observations with $N=4\,000$ sampled at equidistant times with sampling interval $\Delta t=2.5\cdot 10^{-2}$. Figure~\ref{fig:remove}(a) shows the increase of the normalised cost function $C(\Xi)$ upon removal of members of the library $\boldsymbol{\Theta}(\X)$ for fixed delay time $\tau_s$. The normalization is with respect to  $C(0)$, the value when all library functions are removed, i.e. the error encountered for a rough model with a constant solution $x(t)=x(0)$. 
The member of the library $\mathcal{R}$ to be removed at each iteration is chosen as the one leading to the least increase of the normalized cost function. This iterative process terminates with the remaining terms as output, just before the normalized cost function increases by more than 10$\%$.  We present results for the delay time $\tau_s=7$ (blue curve with open circles), corresponding to the true delay time, and for a non-optimal delay time $\tau_s=6$ (red curve with diamonds). We indicate for both delay times the library functions which are removed at each iteration. For the correct delay time $\tau_s=7$, as expected, we observe a jump in the cost function when one of the terms is being removed which appears in the DDE (\ref{e.ENSO}) (i.e. $x(t)$, $x(t)^3$, $x(t-\tau)$). For the non-optimal delay time $\tau_s=6$ we also see, as expected, a jump but the selected terms $x(t)$, $x(t)^3$, $x(t-\tau)^3$ do not correspond to the actual terms appearing in (\ref{e.ENSO}). We also observe a significantly lower value of the cost function for the (optimal) delay time $\tau_s=7$ compared to the non-optimal delay time $\tau_s=6$ at iteration numbers for which none of the selected monomials have been removed. In Figure~\ref{fig:L2error} (open circles, online blue), we show how the optimal delay time $\tau^\star$ is determined by inspecting the reconstruction error $\mathcal{E}(\tau_s)$ (cf. (\ref{e.ell2})). The reconstruction error has a clear minimum at $\tau^\star=7$.
%
In the more challenging case when only ($N=4\,000$) noise-less observations $x(t_n)$ are available and the derivative matrix $\dot \X$ has to be estimated in a post-processing step using the polynomial smoothing described in Remark~\ref{rem:noise}. We perform the polynomial regression with $r=25$ with $\delta = r\Delta t=0.625$. The SINDy algorithm recovers the coefficients and delay times close to the true values as seen in row (b) of Table~\ref{tab:para}.\\

We now test the method in the difficult case of short noise-contaminated data with $N=200$. The variable $x(t)$ is sampled at observation intervals of $\Delta t=0.25$ and are contaminated with observational noise with $\gamma=0.02$. To estimate the derivatives from the data polynomial regression is employed with $r=5$ and $\delta = r\Delta t=1.25$.
Note the smaller value of $r$ compared to the noiseless case considered above accounting for the ten-times larger sampling time used here. Figures~\ref{fig:remove}(b) and \ref{fig:L2error} (diamonds, online red) show that, remarkably, SINDy identifies the correct members of the library and provides an excellent estimation of the delay time with $\tau^\star = 7$. The estimated parameters for the SINDy model (\ref{e.SINDy}) are reported in row (c) of Table~\ref{tab:para}, and, unsurprisingly, more strongly deviate from the true values with reconstruction error $\mathcal{E}(\tau^\star=7)$ of 2.56\%. If a longer time series with $N=4\,000$ is used to train the SINDy model, the error is reduced to 1.94\% (Table~\ref{tab:para}, row (d)), indicating that the limiting factor is the noise rather than the length of the time series. Figure~\ref{fig:remove}(b) shows again the normalized cost function $C(\Xi)$ upon removal of members of the library for fixed delay time $\tau_s$. We show results for the true delay time $\tau_s=7$ and a for a non-optimal delay time $\tau_s=6$. The increase of the cost function upon removal of terms appearing in the true model is less pronounced than in the ideal case of noiseless observations (cf. Figure~\ref{fig:remove}(a)). We remark that the value of the cost function for iteration numbers before the removal of the monomials of the DDE (\ref{e.ENSO}) is significantly larger than for the ideal noiseless case. 
Figure~\ref{fig:L2error} shows the reconstruction error $\mathcal{E}(\tau_m)$ as a function of the delay time $\tau_s$. As in the noiseless case a clear minimum is observed at $\tau_s=7$ corresponding to the delay of the true model. Near the minimum the reconstruction error ${\mathcal{E}(\tau_s)}$ is continuous. For delay times far away from the minimum the reconstruction error experiences discrete jumps, which are caused by the discrete removal of library terms for those values. The perfect accuracy of the estimation of the delay with $\tau^\star=7$ is due to the coarse sampling time $\Delta t=0.25$ implying that the next closest values of delays used for the optimization are $\tau_s=6.75$ and $\tau_s=7.25$, which both lead to a significantly larger value of the reconstruction error $\mathcal{E}(\tau_s)$. We remark that one may use interpolation to provide a more accurate estimate for the delay time at which the minimum of the reconstruction error is attained. The minimum will then not be attained necessarily at a multiple of the sampling time.\\

\begin{figure}[tbp]
\centering
		\begin{subfigure}[t]{0.49\textwidth}
			\centering
			\includegraphics[width=1\linewidth]{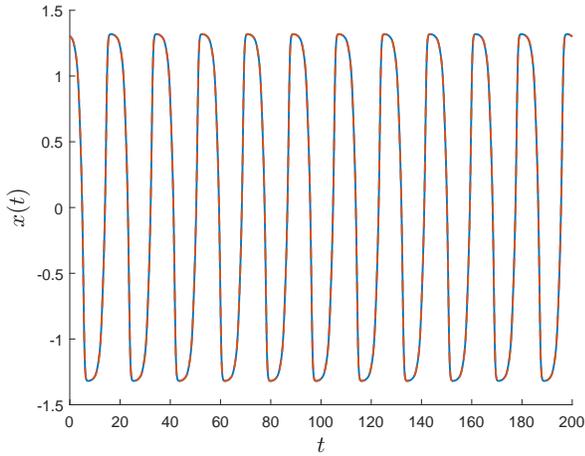}
			\caption{noiseless observations of $x$ and $\dot{x}$\\ \tab $N=4\,000$, optimal $\tau=7$.}
		\end{subfigure}
		\hfill
		\begin{subfigure}[t]{0.49\textwidth}
			\centering
			\includegraphics[width=1\linewidth]{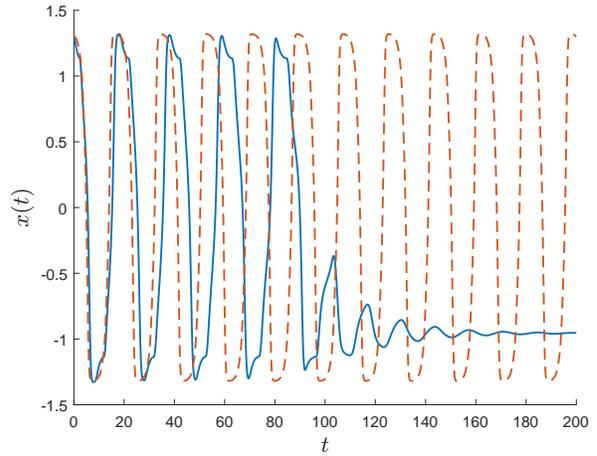}
			\caption{noiseless observations of $x$ and $\dot{x}$\\ $N=4\,000$, non-optimal $\tau=6$.}
		\end{subfigure}\\
		\begin{subfigure}[t]{0.49\textwidth}
			\centering
			\includegraphics[width=1\linewidth]{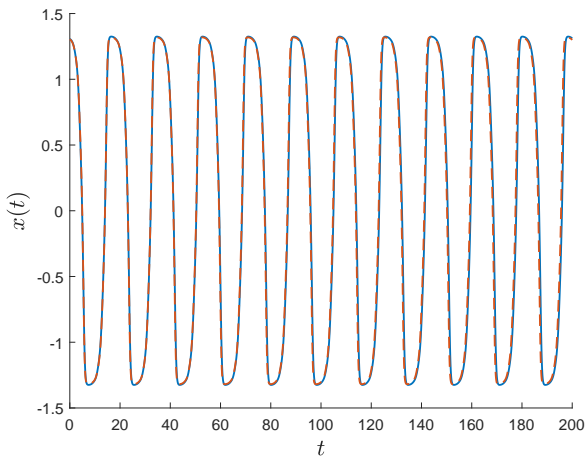}
			\caption{noisy observations of $x$\\ $N=200$, $\gamma=0.02$, optimal $\tau=7$.}
		\end{subfigure}
		\hfill
		\begin{subfigure}[t]{0.49\textwidth}
			\centering
			\includegraphics[width=1\linewidth]{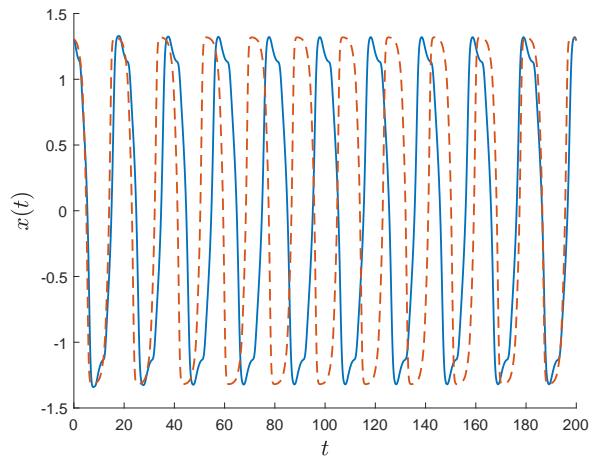}
			\caption{noisy observations of $x$\\ $N=200$, $\gamma=0.02$, non-optimal $\tau=6$.}
		\end{subfigure}
\caption{Comparison of the trajectories obtained from the SINDy model (\ref{e.SINDy}) (continuous curves; online blue) and from the true model (\ref{e.ENSO}) (dashed curves; online red) for optimal and non-optimal delay times, and for noiseless and noisy data points, respectively.}
\label{fig:comparison}
\end{figure}

In Figure~\ref{fig:comparison} we display the trajectories obtained from simulating the estimated SINDy model (\ref{e.SINDy}) and compare them to the trajectory of the (noiseless) true model (\ref{e.ENSO}), both initialized with the same initial condition as the true solution of (\ref{e.ENSO}). Note that the initial value $x(0)=1$ was not part of the (noisy) observations used for training. 
Remarkably, the SINDy algorithm permits to recover the true solution for the observed time window even for the noise-contaminated case with trajectories which are hardly discernible with the bare eye. For longer times we will, however, observe increasing phase errors for the noisy case which does not recover the true coefficients exactly (not shown for brevity). The same SINDy model run with a close but non-optimal delay time of $\tau=6$, however, leads to strong phase and amplitude errors of the SINDy-DDE model as seen in Figures~\ref{fig:comparison}(b) and~\ref{fig:comparison}(d). 

The results presented for the simple one-dimensional toy model (\ref{e.ENSO}) suggest that the SINDy-delay Algorithm~\ref{algo1} described in Section~\ref{sec:SINDy} is able to recover the dynamics of a DDE for relatively short data contaminated by moderate measurement noise at least on the time-scale covered by the observations. Indeed, the reconstruction error is only $2.56$ \% (see Table \ref{tab:para}) with only $N=200$ data measurements, compared to standard applications of SINDy with for instance with $N=10^5$ data measurements for the Lorenz attractor example in \cite[Appendix 4.2]{BruntonEtAl16}. In the next Section we show how to use the SINDy-delay method to uncover the dynamics from a set of biological experiments where the underlying dynamical system is not known.

\section{Application to experimental data of gene expressions in {\em{Pseudomonas aeruginosa}}}
\label{sec:bio}


\begin{figure}[tbp]
\centering
\includegraphics[width=1\columnwidth]{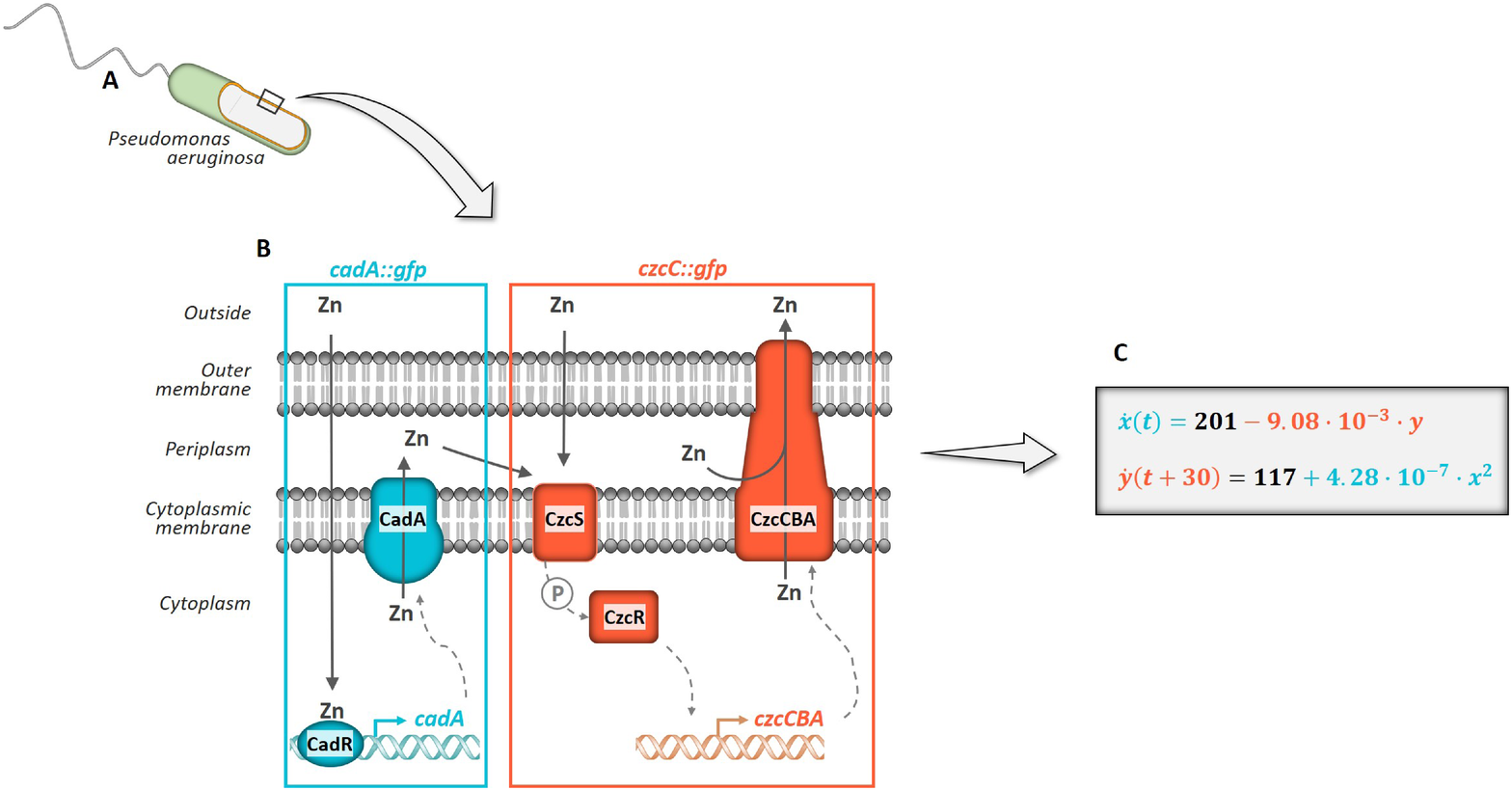}
\caption{ 
Sketch of the biological model of the two compartment system (TCS) for the regulation of zinc in \textit{Pseudomonas aeruginosa}. \textbf{A)} Representation of the bacterium \textit{Pseudomonas aeruginosa}. The square represents the location of the two membranes in which the transport systems visible in \textbf{B} are integrated in. \textbf{B)} Schematic representation of the two-steps
dynamical response of the proteins CadA (blue) and CzcCBA (red) after zinc induction, adapted from \cite{ducret_czccba_2020}. As soon as the metal enters the cell, CadA is rapidly expressed by CadR, leading in a second
phase to the induction of CzcCBA via the CzcRS TCS. \textbf{C)} The delay differential equation describing the dynamics of the Cad (blue color) and Czc (red color) systems after addition of 2 mM Zn obtained by the SINDy-delay method. 
}
\label{fig:biomodel}
\end{figure}

Zinc is an essential element in most living organisms and its proper dosage is vital for their survival. In bacteria, zinc is typically bound to proteins and is responsible for both structural and functional roles of those proteins \cite{blencowe2003}. Too small zinc concentrations impede on the biological functioning of these proteins. Equally, if zinc is present in excess, it becomes toxic, mainly by nonspecific
bindings compromising the cellular integrity \cite{blencowe2014}. 
Therefore, intracellular zinc concentration must be tightly regulated. 
This balance of cellular concentration (also called homeostasis) is finely controlled by zinc import and export systems and their regulators.

Several strategies have evolved in \textit{P. aeruginosa} to mitigate against strong fluctuations of environmental zinc concentrations. In particular, numerous systems composed of transmembrane complexes act to maintain zinc homeostasis \cite{pederick2015,ducret2021}. 
Like all Gram negative bacteria, \textit{P.aeruginosa} possesses a double membrane separated by a particular compartment, the periplasm as illustrated in Figure~\ref{fig:biomodel}. Several complexes are involved in the uptake of zinc in two stages: the first one allows transport of zinc from the outside into the periplasm, the second allows for transport from the periplasm into the cytoplasm. In presence of zinc excess, the associated import transporters are repressed, giving way to a reversed export systems. The most effective transporter is the efflux pump CzcCBA, which expels metal from the periplasm or the cytoplasm directly out of the cell (cf. Figure~\ref{fig:biomodel}) \cite{nies1989,goldberg1999}. The P-type ATPase CadA on the other hand expels zinc from the cytoplasm to the periplasm \cite{lee2001}. (We follow here the standard convention that proteins have names starting with a capital letter whereas their associated genes have names all in small caps and are written in italics). Other export systems have been described in this bacterium, such as CzcD or YiiP, but do not appear to play a major role in zinc resistance \cite{salusso2017,ducret_czccba_2020}.

The expression of the proteins CadA is regulated by CadR that belongs to a family of transcriptional regulators known to be constitutively located on the promoter sequences of their target genes \cite{brown2003}. This configuration provides a fast response as follows: when the cytoplasmic zinc concentration reaches a critical value, CadR binds the metal and immediately induces \textit{cadA} transcription 
\cite{ducret_czccba_2020}. 
The threshold of zinc concentration for the activation of this system depends directly on the zinc affinity of CadR. This threshold has not yet been determined in \textit{P. aeruginosa} in the literature, but it may be estimated about $3 \cdot 10^{-12}$ M, as observed in other bacteria for ZntR, a CadR homolog \cite{Osman2015}. 
Conversely, the efflux pump CzcCBA is activated by the CzcRS TCS, where in presence of high periplasmic concentration of zinc, the CzcS sensor activates the CzcR regulator which in turn binds the DNA, promoting the activation of its own transcription and the \textit{czcCBA} efflux pump, but also represses \textit{oprD} porin transcription~\cite{ducret_czccba_2020}. 
OprD is the entry route for carbapenem antibiotics. Therefore, in presence of zinc, CzcR render the bacterium resistant to both metal and antibiotics. Interestingly, the CadA P-type ATPase appeared to be a key component for a full and timely induction of CzcCBA, suggesting a hierarchical expression
in zinc export systems \cite{ducret_czccba_2020} as shown schematically in Figure~\ref{fig:biomodel}. In a zinc deficient medium, all import systems are expressed.
Consequently, zinc accumulates rapidly in the cytoplasm during a metal boost. This results in the closure of the uptake machineries and at the same time the fast induction of CadA, which begins to expel zinc from the cytoplasm to the periplasm, leading subsequently to the activation of the CzcRS TCS and therefore of CzcCBA. This subsequently promotes a strong expulsion of zinc, which in turn decreases CadR activity and hence CadA expression. To better characterize and model this regulatory system, we seek a simplified two-dimensional differential equation system, describing the dynamical induction of the two agents CadA and CzcCBA. The following subsection describes the experimental set-up employed to obtain measurements for CadA and CzcCBA.

\subsection{Experimental design and results}
\label{subsec:bio_mod}
We used the transcriptional fusions \textit{cadA::gfp} and \textit{czcCBA::gfp} described in \cite{ducret_czccba_2020}. 
This method has the advantage of closely reflecting the expression of the gene of interest and naturally yields time series of experimental data. 
To do so the green fluorescent protein (GFP) were fused with the regulatory sequences of \textit{cadA} or \textit{czcC} genes, respectively. To investigate the interaction between the proteins CadA and CzcCBA, we consider a wild type (wt) strain of \textit{P. aeruginosa} as well as mutants in which either CadA is not expressed ($\Delta$\textit{cadA}) or CzcCBA is not expressed ($\Delta czcA$). Strains were independently grown in a zinc deficient M-LB medium as described in \cite{ducret_czccba_2020}, for 2 hours 30 minutes before the addition of different concentrations of zinc (in the form of ZnCl\textsubscript{2}). We perform experiments for various zinc concentrations with 0.5, 1, 1.25, 1.5, 1.75, 2, 2.25, 2.5 mM, 
a range of zinc concentrations for which the considered systems are fully induced. The fluorescence of \textit{cadA::gfp} was measured for the wt and the $\Delta$\textit{czcA} strains. Similarly, the fluorescence of \textit{CzcCBA::gfp} was measured in the wt and the $\Delta$\textit{cadA} strains. 
Fluorescent values were monitored every 5 minutes for 160 minutes and normalized by the optical density at 600nm (OD600, a standard methodology which permits to estimate the bacterial concentrations). This amounts to a short time series of 33 measurements per experiment. Each experiment is conducted three times and we report on the averages over those three experiments. In the following time $t = 0$ corresponds to the moment of the metal addition. For ease of exposition, fluorescence measurement are shifted to start with a value of $0$ at time $t = 0$. 

Figure~\ref{fig:fusion2mM} shows the fusion measurements for the wild type and two mutants after adding 2mM of ZnCl\textsubscript{2}. In agreement with previous work \cite{ducret_czccba_2020}, in the wt strain (see Figure~\ref{fig:fusion2mM}a), the CadA induction drops when CzcCBA begins to be expressed, i.e. several minutes after the addition of zinc.
However in the $\Delta czcA$ mutant we observe a continuous induction of
CadA during the time of the experiments (see Figure~\ref{fig:fusion2mM}b). The fusion results also reveal a later induction of CzcCBA in the $\Delta cadA$ strain compared to the wt strain (see Figure~\ref{fig:fusion2mM}c).

\subsection{SINDy-delay method to uncover CadA and CzcCBA system dynamics}

\begin{figure}[tbp]

		\begin{subfigure}[b]{0.3\textwidth}
			\raggedleft
			\includegraphics[height=4.5cm]{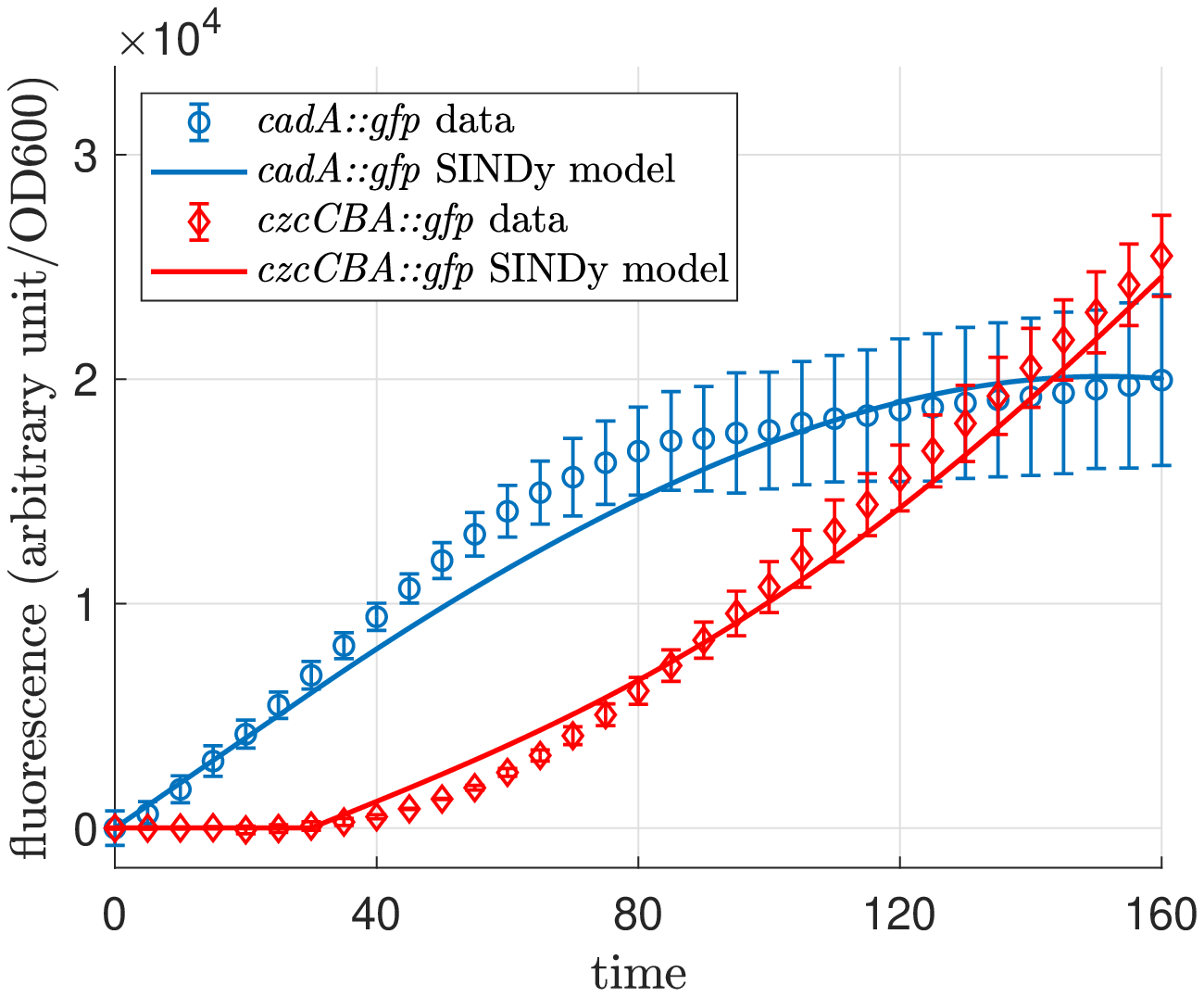}
			\caption{wt}
		\end{subfigure}\quad\quad\quad
		\begin{subfigure}[b]{0.3\textwidth}
			\centering
			\includegraphics[height=4.3cm]{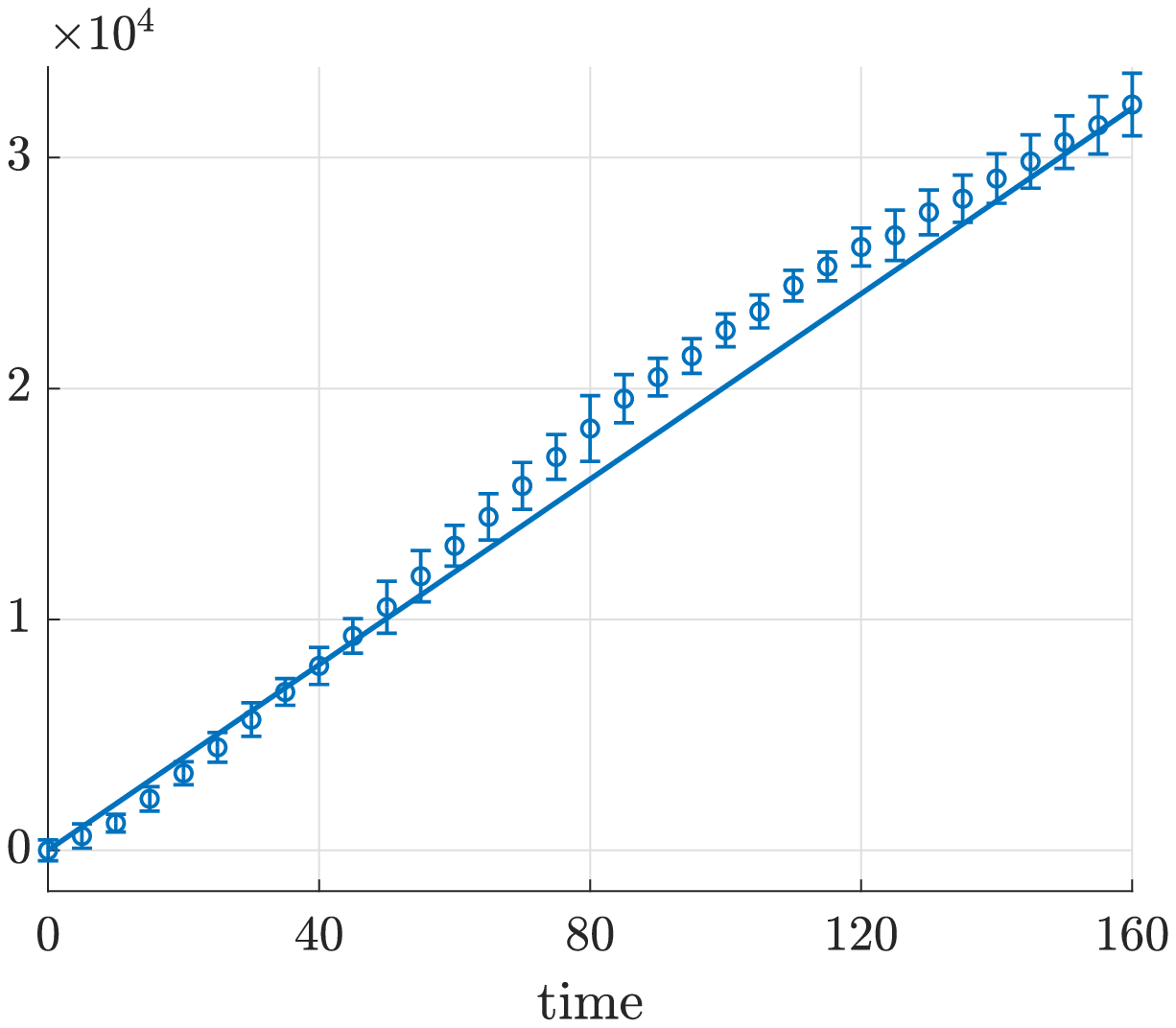}
			\caption{$\Delta$\textit{czcA}}
		\end{subfigure}\quad
		\begin{subfigure}[b]{0.3\textwidth}
			\raggedright
			\includegraphics[height=4.3cm]{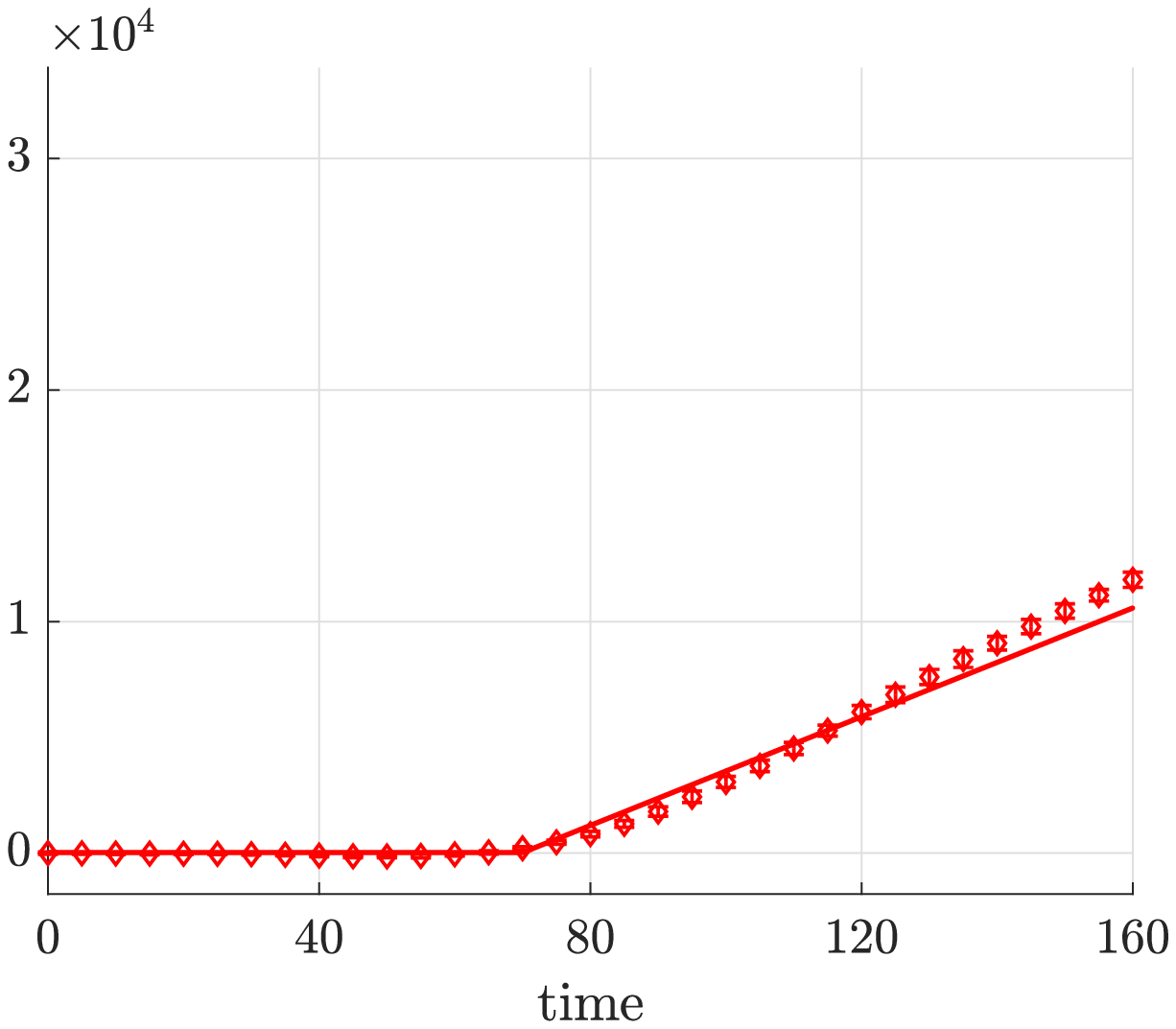}
			\caption{$\Delta$\textit{cadA}}
		\end{subfigure}

\caption{
Fluorescence measurements after addition of 2 mM ZnCl\textsubscript{2} compared to the corresponding mathematical delay differential equation (DDE) model (solid lines). The fluorescence intensity over time is shown for the wt, $\Delta cadA$ and $\Delta czcA$ strains, containing the \textit{cadA::gfp} (open circles; online blue) or \textit{czcA::gfp} (diamonds; online red) fusions. The values are normalized by the optical density (OD600). Standard deviations of three independent measurements are shown. The mathematical solutions, according to the SINDy selected model, are shown in solid lines.
}
\label{fig:fusion2mM}
\end{figure}

The dynamics and induction intensity of CadA and CzcCBA systems depend on several factors, including intracellular (periplasmic and/or cytoplasmic) concentration of zinc, as well as on the response velocity and the metal sensitivity of their respective regulators. Moreover, experimental data obtained 
from transcriptional fusions are only proxies depending on GFP synthesis and its stability.
For simplicity, we ignore these complex interactions and instead consider only two “boxes”, one signifying all the variables involved in CadA expression (blue box in Figure~\ref{fig:biomodel}) and one signifying those responsible of CzcCBA expression (red box in Figure~\ref{fig:biomodel}). This simplification implies a mathematical model with only CadA and CzcCBA expressions as dependent variables. We assume that the zinc concentration remains constant during the induction experiment. 

For the wild type bacteria (wt) we seek a model of the form
\begin{align} 
\label{e.frawark1a}
\dot x_{wt}(t) &= f(x_{wt},y_{wt}), \qquad \dot y_{wt}(t+\tau_{wt}) = g(x_{wt},y_{wt}),
\end{align}
where $x(t)$ represents the fluorescence from \textit{cadA::gfp} while  $y(t)$ represents \textit{czcA::gfp}. This form is motivated by the experimental data shown in Figure~\ref{fig:fusion2mM} where \textit{cadA::gfp} experiences significant changes within the first minutes whereas \textit{czcA::gfp} remains nearly constant for a significant time suggesting a delayed dynamics. For the mutant $\Delta\textit{czcA}$, which lacks expression of $\textit{czcA}$ the dynamics is obtained by setting $y=0$ in the above model for the wild type. We obtain
\begin{align}
\dot x_{\Delta\textit{cz}}(t) &= f(x_{\Delta\textit{cz}},0) .
\label{e.frawark1b}
\end{align}
Similarly, for the mutant $\Delta\textit{cadA}$, which lacks expression of $\textit{cadA}$ the dynamics is obtained by setting $x=0$ in the above model for the wild type, and we obtain
\begin{align}
\dot y_{\Delta\textit{ca}}(t+\tau_{\Delta\textit{ca}}) &= g(0,y_{\Delta\textit{ca}}).
\label{e.frawark1c}
\end{align}
We allowed here for a delay time $\tau_{\Delta\textit{ca}} \neq \tau_{wt}$ accounting for the possibility that the delay may depend on the presence of the various agents present in the regulatory process.
We also assume that $x(t)=y(t)=0$ for $t<0$ which corresponds to the natural assumption that neither CadA ($x$) nor CzcCBA ($y$) are produced when no zinc has been added yet into the growing medium, which could activate their expression (see Figure~\ref{fig:biomodel}).

%

To determine the model \eqref{e.frawark1a}-\eqref{e.frawark1c}, we apply the SINDy-delay method presented in Sections~\ref{sec:SINDy} for the fluorescence measurements of the expression kinetics experiments described in Section~\ref{subsec:bio_mod}. To estimate the functions $f$ and $g$ in \eqref{e.frawark1a} as well as the delay times from the experimental data, we consider a library consisting of all monomials up to cubic order to approximate \eqref{e.frawark1a}-\eqref{e.frawark1c}. 
To search for a parsimonious model only few terms selected from the library, we apply the sparsity constraints as detailed in Section \ref{sec:SINDy} and search for the delay times $\tau_{\Delta\textit{ca}}$ and $\tau_{wt}$ in the set $[0,5,\ldots,160]$ in units of minutes: we remove the less meaningful terms of the library and stop the process when the minimum of the normalized cost function $C(\Xi)$/$C(0)$, which refers to equation (\ref{e.C}), increases by more than 10 percents. Optimal delays time are found by minimizing the function $\mathcal{E}(\tau_{wt},\tau_{\Delta ca})$ reconstruction error corresponding to (\ref{e.ell2}). 

\begin{figure}[tbp]
		\begin{subfigure}[b]{0.5\textwidth}
			\raggedleft
			\includegraphics[width=1\linewidth]{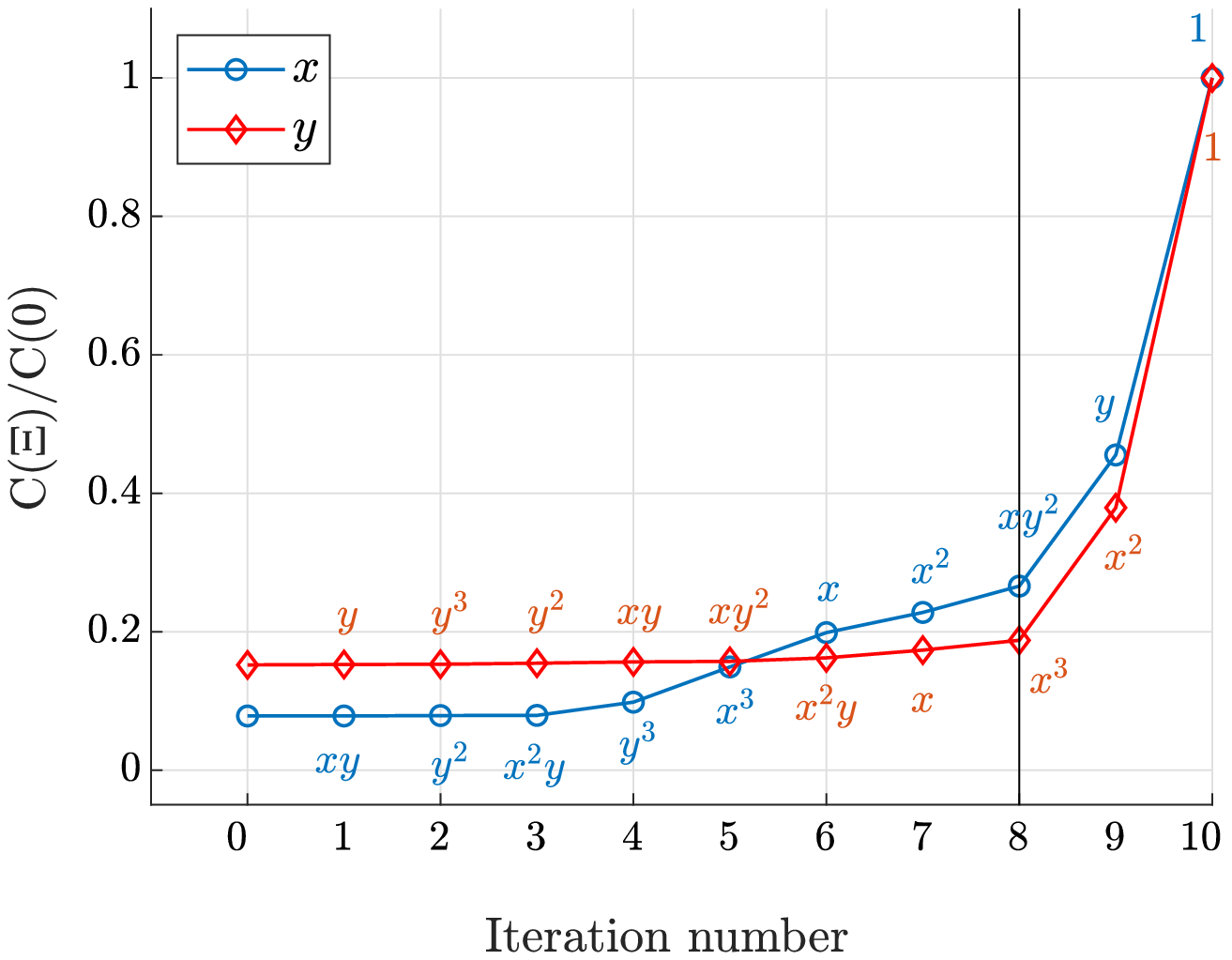}
			\caption{}
		\end{subfigure}
		\begin{subfigure}[b]{0.5\textwidth}
			\includegraphics[width=1\linewidth]{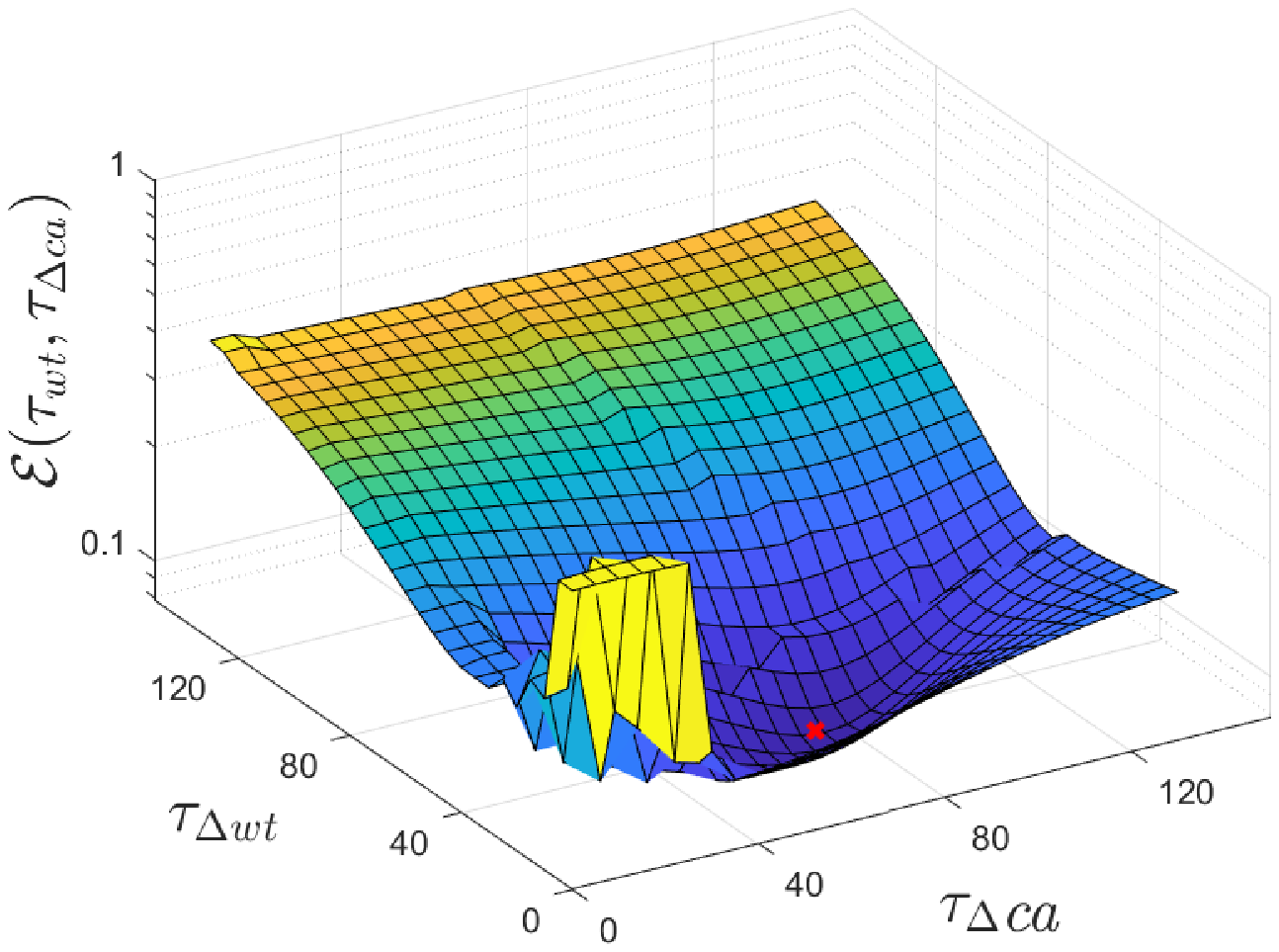}
			\caption{}
		\end{subfigure}
\caption{(a) Cost function $C(\Xi)$ (normalized by $C(0)$) against removed monomials for fixed optimal delays, $\tau_{wt}=30$ (open circles; online blue) and $\tau_{\Delta ca}=70$ (diamonds; online red). The vertical black line indicates the iteration number where the process is stopped. At iteration 8 the cost function has increased more than 10\% for both components. (b) Reconstruction error $\mathcal{E}(\tau_{wt},\tau_{\Delta ca})$ showing a minimum value equal to $7.8\cdot10^{-2}$ for $\tau_{wt}=30$ and $\tau_{\Delta ca}=70$, indicated by a red cross.}
\label{fig:2mMprocess}
\end{figure}

This process is applied for all experiments with the various zinc concentrations. In Figure~\ref{fig:2mMprocess} we show results for the 2 mM induction of zinc. 
Figure~\ref{fig:2mMprocess}(a) shows the increase of the normalized cost function upon removal of both $x$ and $y$ components. 
Figure~\ref{fig:2mMprocess}(b) shows the reconstruction error $\mathcal{E}(\tau_s)$ with a minimum error of $7.8\%$ for $\tau_{wt}=30$ and $\tau_{\Delta\textit{ca}}=70$ minutes. 
In particular, we obtain from the SINDy-delay methodology the following delay differential equation model for a concentration of 2 mM of zinc, 
\begin{align}\centering
\dot x_{wt}(t) &= 201-9.08\cdot10^{-3}\, y_{wt}, \label{eq:supermodelx}\\
\dot y_{wt}(t+\tau_{wt}) &= 117+4.28\cdot 10^{-7}\, x_{wt}^{2}
\label{eq:supermodely}
\end{align}
and 
\begin{align}\centering
\dot x_{\Delta cz} &=201, \nonumber\\
\dot y_{\Delta ca}(t+\tau_{\Delta ca}) &=117
\label{eq:supermodel3}
\end{align}
with $\tau_{wt}=30, \tau_{\Delta ca}=70$.
In Figure~\ref{fig:fusion2mM}, solutions of the DDE model \eqref{eq:supermodelx}-\eqref{eq:supermodely} and of \eqref{eq:supermodel3} are plotted and compared with experimental data for 2 mM ZnCl\textsubscript{2}, which shows a high degree of similarity with a reconstruction error of 7.85\%.  
The complete results for all zinc concentrations tested (from 0.5 to 2.5 mM) are shown in Table \ref{tab:parabio}. We remark that the coefficient of the linear and quadratic terms in \eqref{eq:supermodelx},\eqref{eq:supermodely}, are of the order of $10^{-2}$ and $10^{-7}$, respectively; although their coefficients are small, their presence is crucial. Such small coefficients are hard to detect when employing standard thresholding procedures. This illustrates the advantage of our method based to promote sparsity outlined in Remark~\ref{rem:sparsity}.

\begin{table}[btp]
\begin{center}
\renewcommand{\arraystretch}{1.3}
\scalebox{0.6}{
\begin{tabular}{|c|cc|ccccccc|ccccc|c|}
\hline
Zn &\multicolumn{2}{|c|}{Delays} & \multicolumn{7}{|c|}{Coefficients for function $x(t)$  (\textit{cadA::gfp})} & \multicolumn{5}{|c|}{Coefficients for function $y(t)$  (\textit{czcA::gfp})}&Error \\
\hline
$[mM]$ & $\tau_{wt}$& $\tau_{\Delta\textit{ca}}$ & $1$ & $x$  &$y$&$x^2$&$xy$ & $y^2$ & $xy^2$& $1$& $x$&$y$&$x^2$ & $xy$ & $\mathcal{E}(\tau_{wt},\tau_{\Delta ca})$ \\
\hline
0.5 &20  &25 &370 &-3.51$\cdot10^{-2}$  & -3.70$\cdot10^{-2}$  &1.14$\cdot10^{-6}$	&0	& 1.69$\cdot10^{-6}$&0&118&1.61$\cdot10^{-2}$&0&0&-1.53$\cdot10^{-6}$&0.0328\\
1&	25 & 35	& 317&	0&	0&-2.16$\cdot10^{-7}$&	-3.55$\cdot10^{-6}$&0&1.14$\cdot10^{-10}$	&180&5.30$\cdot10^{-3}$	&-6.51$\cdot10^{-3}$&0&0&0.0711\\
1.25&30&45&316&0&0&-1.74$\cdot10^{-7}$&-2.56$\cdot10^{-6}$&0&7.90$\cdot10^{-11}$&162&0&-5.18$\cdot10^{-3}$&3.36$\cdot10^{-7}$&0&0.0666\\
\hdashline
{1.5}&	25&	{55}&	216&0&	-1.13$\cdot10^{-2}$&0&0&0&0&122&0&0&2.35$\cdot10^{-7}$&0&	0.11\\
{1.75}&	30&	{65}&	209&0&	-1.04$\cdot10^{-2}$&0&0&0&0&112&0&0&3.46$\cdot10^{-7}$&	0&	0.0865\\
{2}&	30&	{70}&	201&0&-9.08$\cdot10^{-3}$&0&0&0&0&117&0&0&4.28$\cdot10^{-7}$&0&	0.0785\\
{2.25}&	35&	{85}&	200&0&-8.48$\cdot10^{-3}$&0&0&0&0&134&0&0&	4.10$\cdot10^{-7}$&0&0.0730\\
\hdashline 
2.5&	45&	95&	200&0&-6.54$\cdot10^{-3}$&0&0&0&0&129&	1.04$\cdot10^{-2}$&0&0&0&	0.0613\\
\hline
\end{tabular}}
\end{center}
\caption{Results of the SINDy-delay method for the various zinc concentrations. Terms from the library function which were not selected for any zinc concentration are not represented. 
}
\label{tab:parabio}
\end{table}%
\renewcommand{\arraystretch}{1}%

\paragraph{DDE model accuracy and consistency} 
We observe in Table~\ref{tab:parabio} that for all zinc concentrations the SINDy-delay method yields DDE models with reconstruction errors smaller than 11\%. The SINDy model matches the experimental data very well and is biologically consistent for all ZnCl\textsubscript{2} concentrations. This is notable given the very short length of the experimental time series with $N=33$ data. 
Remarkably, for moderate ZnCl\textsubscript{2} concentrations between 1.5 mM and 2.25 mM (emphasized in dashed lines in Table~\ref{tab:parabio}), a unified SINDy DDE model arises
which benefits from the sparsity feature, with only the terms $1$,$y$ and $x^2$ selected, allowing for a biologically consistent interpretation of the terms. 
Importantly, the signs of the associated coefficients are consistent with the biological model: the coefficient associated with the linear term in $y$ in (\ref{eq:supermodelx}), which describes the influence of $y$ (CzcCBA) on $x$ (CadA), is negative in agreement with the biological model where CzcCBA represses CadA. Similarly, the coefficient associated with the $x^2$ term in (\ref{eq:supermodely}), which describes the influence of $x$ (CadA) on $y$ (CzcCBA), is positive in agreement with the biological model where CadA accelerates the expression of CzcCBA \eqref{e.frawark1a}.

For low ZnCl\textsubscript{2} concentrations smaller than 1.25 mM the SINDy models are not as sparse, involving more terms than for the moderate concentrations (for instance, up to five functions $1,x,y,x^2,y^2$ for $x$ (CadA) at 0.5 mM of zinc), while for the highest considered ZnCl\textsubscript{2} concentration, the different term $x$ is selected in place of $x^2$.
We also remark that the SINDy model is likely to model the response of \textit{P. aeruginosa} to a boost in zinc only for the time duration of the experiment. Indeed, the SINDy models in Table~\ref{tab:parabio} exhibit unphysical negative CadA and CzcA concentrations for all considered ZnCl\textsubscript{2} concentrations if a time larger than 800 minutes would be considered (not shown here) in place of the time of 160 minutes considered in the experiments. 
This suggests that the simplified two box model may be insufficient to capture the impact of the induced stress for longer times, and additional components or mechanisms need to be included in the modelling.

\begin{figure}[tbp]
\centering
			\includegraphics[width=0.5\linewidth]{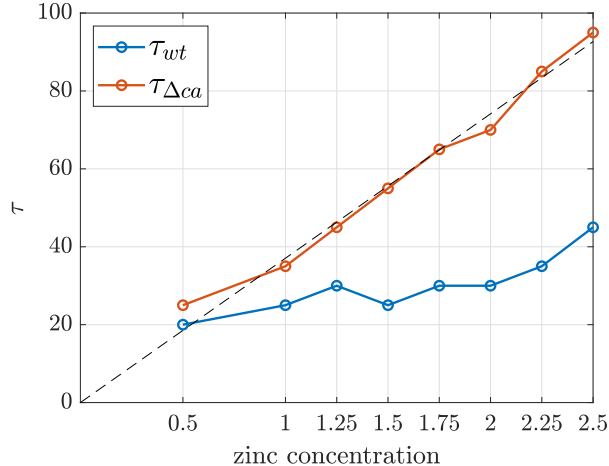}
\caption{Estimated delay times $\tau_{wt}, \tau_{\Delta ca}$ (in minutes) as a function of the zinc concentration (in mM). Remarkably, we observe a linear increase as a function of the zinc concentration of the delay time $\tau_{\Delta ca} = \alpha \cdot [\mathrm{ZnCl_\textsubscript{2}}]$ with $\alpha = 37.1$ min mM\textsuperscript{-1} (linear regression in dashed line with slope $37.1$).} 
\label{fig:delaycomp}
\end{figure}

\paragraph{CadA is essential for maintaining a rapid expression of CzcCBA}
Consider the range of zinc concentrations from 1.5 to 2.25 mM as emphasised with dashed lines in Table~\ref{tab:parabio}. A remarkable observation is that the coefficients computed from the SINDy-delay method are only weakly sensitive to the applied zinc concentration, with the exception of the delay time $\tau_{\Delta ca}$, which increases linearly with the zinc concentration, as shown in Figure \ref{fig:delaycomp}.
This linear increase of the delay time $\tau_{\Delta ca}$ in the absence of CadA suggests that the protein CadA is particularly necessary for a rapid zinc response 
and suggests that the positive effect of CadA on the efflux pump is all the more important as the zinc concentration is high. 
The OD600 measurement allows the counting of cells independently of whether they are alive or dead. Thanks to colony counting and quantification of cell viability at concentrations
of 1.25 mM and 2 mM ZnCl\textsubscript{2} after 160 min of incubation (not displayed here for brevity), we observed that the same number of living cells 
are detected, and hence this difference in delay under different zinc concentrations cannot be attributed to a differential mortality between the $\Delta cadA$ and the wt strains.
Biologically, this could reflect a reasonable mechanism whereby the bacterium wants to react as quickly as possible to a stress regardless of its intensity.

\section{Conclusion}
\label{sec:conclusion}

In this paper, we extended the SINDy methodology introduced in \cite{BruntonEtAl16} to the case of delay differential equations with a focus on short and noise-contaminated data. 
To construct the temporal derivatives from noisy measurements we employed a simple denoising procedure based on polynomial regression (Remark~\ref{rem:noise}). We further introduced a stopping criterion to promote sparsity which avoids having to introduce sensitive threshold parameters (Remark~\ref{rem:sparsity}). To estimate the temporal delay we applied a bilevel optimization whereby first standard SINDy method is applied for a range of fixed delay times, and then subsequently the optimal delay time is determined by the delay time yielding the minimal reconstruction error. We showed that our method is able to reliably uncover the DDE from noisy data obtained from a known toy model.\\

Applying the SINDy-delay methodology to model the dynamics of the {\em Pseudomonas aeruginosa} zinc response from a limited amount of  measurement 
highlighted the subtle interactions between the Cad and Czc regulatory systems. In particular, the SINDy DDE model revealed the importance of CadA on CzcCBA induction for minimizing the time required for the bacterium to respond effectively to a sudden zinc excess. 
The compatibility between the results of the SINDy DDE models and the biological data 
supports the hypothesis that the dynamical mechanism of resistance to moderate boosts of zinc can be explained by the interaction of only two systems, namely CadA and CzcCBA. Our results motivate further investigations of this dynamics. The present work was performed over 160 minutes after the metal induction and illustrates only the initial establishment of resistance. 
Additional experimental data on longer times, which require continuous cultures in a chemostat and a more sensitive method to monitor the \textit{cadA} and \textit{czcCBA} transcriptional expressions would make it possible to compare these mathematical predictions with the biological situation.



\paragraph{Acknowledgments}
This work was partially supported by the UNIGE-USyd strategic Partnership  Collaboration Awards (PCA) 2019-2024, the Swiss National Science
Foundation, project No. 31003A\_179336 for K.P and projects No. 200020\_184614 and No. 200020\_192129 for G.V. 
G.A.G. acknowledges funding from The Australian Research Council, grant DP220100931.



\bibliographystyle{abbrv}
\bibliography{SINDy} 

\end{document}